\documentclass[12pt]{article}
\usepackage{a4}
\usepackage{amsthm}
\usepackage{amsfonts}
\usepackage{amssymb}
\usepackage{cite}
\usepackage{epsfig}
\newtheorem{conjecture}{Conjecture}
\newtheorem{theorem}{Theorem}
\newtheorem{lemma}[theorem]{Lemma}

\begin{document}
\title{Steinberg's Conjecture is false}
\author{Vincent Cohen-Addad\thanks{D\'epartement d'Informatique, UMR CNRS 8548, \'Ecole Normale Sup\'erieure, Paris, France. E-mail: {\tt \{vcohen,zhentao\}@di.ens.fr}.}\and
\newcounter{lth}
\setcounter{lth}{5}
        Michael Hebdige\thanks{Department of Computer Science, University of Warwick, Coventry CV4 7AL, UK. E-mail: {\tt m.hebdige@warwick.ac.uk}. The work of this author was supported by the Leverhulme Trust 2014 Philip Leverhulme Prize of the third author.}\hskip 1ex $^\fnsymbol{lth}$\and 
        Daniel Kr\'al'\thanks{Mathematics Institute, DIMAP and Department of Computer Science, University of Warwick, Coventry CV4 7AL, UK. E-mail: {\tt d.kral@warwick.ac.uk}. The work of this author was supported by the Leverhulme Trust 2014 Philip Leverhulme Prize.}\hskip 1ex $^\fnsymbol{lth}$\and
\setcounter{lth}{1}
        Zhentao Li$^\fnsymbol{lth}$\and
        Esteban Salgado\thanks{\'Ecole Polytechnique, Paris, France, and Pontificia Universidad Cat\'olica de Chile, Santiago, Chile. E-mail: {\tt easalgadovalenzuela@gmail.com}.}\hskip 1ex \thanks{This research was done during the stay of the second, third and fifth of the authors at \'Ecole Normale Sup\'erieure in Paris.}}
\date{}
\maketitle

\begin{abstract}
Steinberg conjectured in 1976 that
every planar graph with no cycles of length four or five is 3-colorable.
We disprove this conjecture.
\end{abstract}

\section{Introduction}

One of major open problems on colorings of planar graphs is Steinberg's Conjecture from 1976.
The conjecture asserts that every planar graph with no cycles of length four or five is 3-colorable.
This problem has been attracting a substantial amount of attention among graph theorists, see~\cite[Section 7]{bib-borodin13}.
It is also one of the six graph theory problems ranked with the 4-star (highest) importance in the Open Problem Garden~\cite{bib-garden}.
As a possible approach towards proving the conjecture,
Erd\H os in 1991 (see~\cite{bib-steinberg93}) suggested to determine the smallest $k$ such that
every planar graph with no cycles of length $4,\ldots,k$ is $3$-colorable.
Borodin et al.~\cite{bib-borodin05+}, improving up on~\cite{bib-abbot91+,bib-sanders95+}, showed that $k\le 7$,
i.e., every planar graph with no cycles of length four, five, six or seven is 3-colorable.
Many other relaxations of the conjecture have been established,
e.g.~\cite{bib-borodin96,bib-borodin03+,bib-borodin09+,bib-xu06}.
On the other hand, the conjecture is known to be false in the list coloring setting~\cite{bib-voigt07}.
We refer the reader for further results on the conjecture and
other problems related to colorings of planar graphs to~\cite{bib-borodin13,bib-jensen95+}.

We disprove Steinberg's Conjecture by constructing a planar graph with no cycles of length four or five that is not 3-colorable.
This also disproves Strong Bordeaux Conjecture from~\cite{bib-borodin03+} and
Novosibirsk 3-Color Conjecture from~\cite{bib-borodin06+}.
\begin{conjecture}[Strong Bordeaux Conjecture]
Every planar graph with no pair of cycles of length three sharing an edge and with no cycle of length five is 3-colorable. 
\end{conjecture}
\begin{conjecture}[Novosibirsk 3-Color Conjecture]
Every planar graph without a cycle of length three sharing an edge with a cycle of length three or five is 3-colorable.
\end{conjecture}
We would like to remark that we have not tried to minimize the number of vertices in our construction;
we are actually aware of the existence of smaller counterexamples (the smallest one with 85 vertices)
but we decided to present the one that is the simplest to analyze among those that we have found.

\section{Counterexample}

Our construction has three steps.
In each of these steps, we construct a particular planar graph with no cycles of length four or five
that has some additional properties.
The next two lemmas present the first two steps of the construction.

\begin{figure}
\begin{center}
\epsfbox{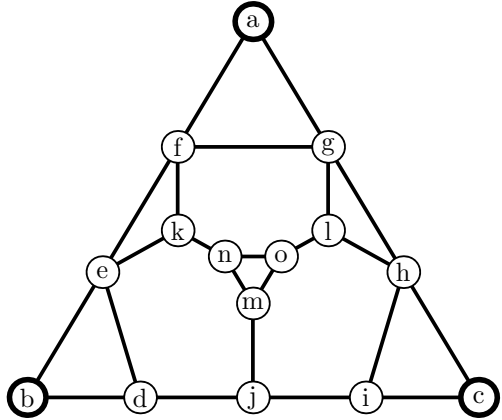}
\end{center}
\caption{The graph $G_1$ analyzed in Lemma~\ref{lm-step1}.}
\label{fig-step1}
\end{figure}

\begin{lemma}
\label{lm-step1}
The graph $G_1$ that is depicted in Figure~\ref{fig-step1} has no cycles of length four or five and
there is no $3$-coloring that would assign the same color to all the three vertices $a$, $b$ and $c$.
In addition, the distance between $a$ and $b$ is three, between $a$ and $c$ is three, and between $b$ and $c$ is four.
\end{lemma}

\begin{proof}
It is easy to check that $G_1$ contains no cycle of length four or five,
the distance between $a$ and $b$ is three, between $a$ and $c$ is three, and between $b$ and $c$ is four.
Suppose that $G_1$ has a $3$-coloring such that all the three vertices $a$, $b$ and $c$ have the same color, say $\gamma$.
The other two colors must be used in an alternating way to color the vertices $d$, $e$, $f$, $g$, $h$ and $i$.
This implies that the vertices $j$, $k$ and $l$ have the color $\gamma$.
Hence, the two colors different from $\gamma$ must be used to color the vertices $m$, $n$ and $o$,
which is impossible since these three vertices form a triangle.
\end{proof}

\begin{figure}
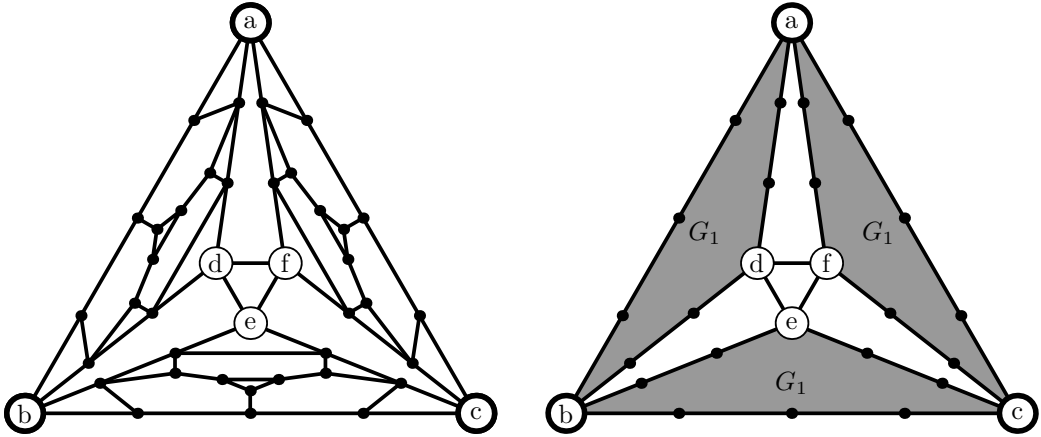

\begin{center}
\epsfbox{steinberg.2}
\hskip 5mm
\epsfbox{steinberg.3}
\end{center}
\caption{The graph $G_2$ analyzed in Lemma~\ref{lm-step2}.
	 The actual graph is drawn on the left;
	 its abstract structure using the graph $G_1$ is given on the right.
         The three contact vertices $a$, $b$ and $c$ are drawn with bold circles.}
\label{fig-step2}
\end{figure}

\begin{lemma}
\label{lm-step2}
The graph $G_2$ that is depicted in Figure~\ref{fig-step2} has no cycles of length four or five and
there is no $3$-coloring that would assign the same color to all the three vertices $a$, $b$ and $c$.
In addition, the distance between any pairs of the vertices $a$, $b$ and $c$ is four.
\end{lemma}

\begin{proof}
Note that the graph $G_2$ consists of three copies of the graph $G_1$ pasted together as depicted in Figure~\ref{fig-step2}.
We start by checking that the distance between the vertices $a$ and $b$ is four (the cases of the other two pairs are symmetric).
Since $a$ and $b$ are joined by a path of length four, their distance is at most four.
Lemma~\ref{lm-step1} implies that there is no path of length less than four inside the copy of $G_1$ containing $a$ and $b$.
On the other hand,
any path passing through more than one copy of $G_1$ must have length at least three inside each traversed copy of $G_1$,
so its total length cannot be less than six.
We conclude that the distance between $a$ and $b$ is four.

We next check that $G_2$ has no cycle of length four or five.
By Lemma~\ref{lm-step1}, no cycle inside one of the copies of $G_1$ has length four or five.
Since the distance between any two of the contact vertices of the copies of $G_1$ (the vertices used to glue the copies together)
is at least three by Lemma~\ref{lm-step1}, the length of any cycle passing through two or three of the copies is at least six.
The only cycle that does not contain an edge inside a copy of $G_1$ is the cycle $def$ and its length is three.
Therefore, $G_2$ has no cycle of length four or five.

It remains to check that $G_2$ has no $3$-coloring such that the vertices $a$, $b$ and $c$ get the same color.
Suppose that $G_2$ has a $3$-coloring such that the vertices $a$, $b$ and $c$ receive the same color, say $\gamma$.
By Lemma~\ref{lm-step1}, the color of any of the vertices $d$, $e$ and $f$ is different from $\gamma$.
However, this is impossible since the vertices $d$, $e$ and $f$ form a triangle and there are only three colors in total.
Hence, $G_2$ has no $3$-coloring such that the vertices $a$, $b$ and $c$ would get the same color.
\end{proof}

We are now ready to provide the counterexample construction.

\begin{theorem}
\label{thm-main}
There exists a planar graph with no cycles of length four or five that is not $3$-colorable.
\end{theorem}

\begin{figure}
\begin{center}
\epsfbox{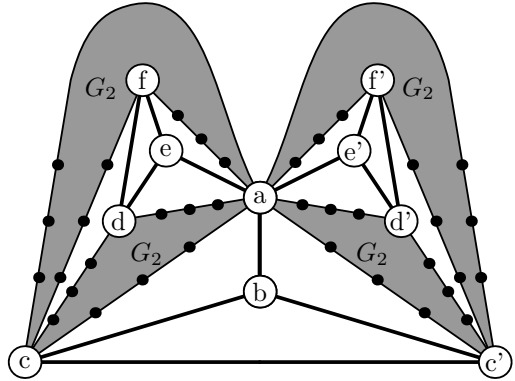}
\end{center}
\caption{The graph $G$ constructed in the proof of Theorem~\ref{thm-main}.
         The gray areas represent the graph from Figure~\ref{fig-step2}.}
\label{fig-main}
\end{figure}

\begin{proof}
Consider the graph $G$ obtained by taking four copies of the graph $G_2$ from Lemma~\ref{lm-step2} and
pasting them together with additional vertices $b$, $e$ and $e'$ in the way depicted in Figure~\ref{fig-main}.
We first argue that $G$ has no cycle of length four or five.
By Lemma~\ref{lm-step2}, there is no cycle of length four or five inside a single copy of $G_2$.
Any cycle containing an edge inside a copy of $G_2$ that is not fully contained inside the copy
must contain at least four edges inside the considered copy of $G_2$,
since the distance between any two of the contact vertices of $G_2$ is four by Lemma~\ref{lm-step2}.
Since such a cycle must contain at least two edges outside the copy of $G_2$,
its length must be at least six.
Hence, there is no cycle  with length four or five that contains an edge inside a copy of $G_2$.
Observe that there are exactly three cycles not containing an edge inside any copy of $G_2$ and
these are triangles $bcc'$, $def$ and $d'e'f'$.
Therefore, $G$ has no cycle of length four or five.

Suppose that $G$ has a $3$-coloring and let $\gamma$ be the color assigned to the vertex $a$.
Since the vertices of the triangle $bcc'$ must be assigned all three colors and the color of $b$ is not $\gamma$,
either $c$ or $c'$ is colored with $\gamma$.
By symmetry, we can assume that the color of $c$ is $\gamma$.
Since the vertices $a$, $c$ and $d$ are contact vertices of a copy of $G_2$,
the color of $d$ is not $\gamma$ by Lemma~\ref{lm-step2}.
Similarly, since the vertices $a$, $c$ and $f$ are contact vertices of a copy of $G_2$,
the color of $f$ is also not $\gamma$.
Hence, none of the vertices $d$, $e$ and $f$ (note the vertex $e$ is adjacent to $a$) has the color $\gamma$,
which is impossible since the vertices $d$, $e$ and $f$ form a triangle.
We conclude that the graph $G$ has no $3$-coloring.
\end{proof}

\section*{Acknowledgement}

The authors would like to thank Zden\v ek Dvo\v r\'ak for his comments and insights
on a computer-assisted approach towards a possible solution of Steinberg's Conjecture.

\end{document}